\documentclass[11pt]{article}

\begin{document}
\noindent
\begin{center}
{\bf A Formula For All $K$-Gonal Numbers That Are Centered $K$-Gonal}
\end{center}

\section{Introduction}
\indent

Figurate numbers probably originated with the Pythagorean society [2, p.58]. These positive integers comprise the total number of dots in various regular polygons (such as triangles, squares, pentagons, hexagons, etc.). Many interesting properties of figurate numbers have been discovered [3]. It is known that the $n^{\mbox{th}}$ triangular number, $t_n$, is given by
$$t_n = \frac{n(n + 1)}{2}$$
and the $n^{\mbox{th}}$ square nunmber is given by
$$s_n = n^2$$
for $n \geq 1$. The first eight triangular numbers are 1, 3, 6, 10, 15, 21, 28, 36. The first six square numbers are 1, 4, 9, 16, 25, 36. So, 1 and 36 are both triangular numbers and square numbers. In fact, it is known that there exist infinitely many numbers that are both triangular and square numbers. In [1], it is shown that the $n^{\mbox{th}}$ triangular square number is
$$\frac{\left[\left(17 + 12\sqrt{2}\right)^n + \left(17 - 12\sqrt{2}\right)^n\right]}{32}.$$
The purpose of this paper is to derive a similar result for a different class of positive integers, namely, $k$-gonal and centered $k$-gonal numbers. For clarity, we introduce a specific example for $k = 3$. 

Geometrically, {\em centered triangular numbers}, $T_n$, for instance, consist of a central dot (or vertex), $T_1 = 1$, with three dots around it and then additional dots in the gaps between adjacent dots (see Figure 1). So, the number of vertices on each edge of each triangle is one larger than the number of vertices on the edge of the preceding smaller triangle. Hence, $T_n$ can be defined recursively by 
\begin{eqnarray*}
T_1 & = & 1\\
T_n & = & T_{n-1} + 3(n - 1)
\end{eqnarray*}
for $ n \geq 2$. The first four centered triangular numbers are 1, 4, 10, 19. From the recursive definition of these numbers, one obtains
$$T_n = \frac{3n^2 - 3n + 2}{2}$$
for $n \geq 1$. Table 1 displays the first six triangular numbers that are centered triangular numbers.

\vspace{1pc}

$$
\begin{array}{c|c}
\mbox{Triangular numbers} & \mbox{Centered triangular numbers}\\
\hline\\
t_1 = 1 & T_1 = 1\\
t_4 = 10 & T_3 = 10\\
t_{16} = 136 & T_{10} = 136\\
t_{61} = 1891 & T_{36} = 1891\\
t_{229} = 26335 & T_{133}= 26335\\
t_{856} = 366796 & T_{495} = 366796
\end{array}
$$
\begin{center}
Table 1
\end{center}

These concepts can be generalized to $k$-gonal [2] and centered $k$-gonal numbers. In general, for a positive integer $k \geq 3$, let $P(n; k)$ be the $n^{\mbox{th}}$ term of the $k^{\mbox{th}}$ polygonal number and $C(n; k)$ be the $n^{\mbox{th}}$ term of the $k^{\mbox{th}}$ centered polygonal number. Then
$$P(n; k) = \frac{(k - 2)n^2 + (4 - k)n}{2}$$
for $n \geq 1$ and
$$C(n; k) = \frac{kn^2 - kn + 2}{2}$$
for $n \geq 0$.
In this paper, we give a specific formula for all $k$-gonal numbers that are also centered $k$-gonal numbers. Our method is based on reducing the problem to solving a Pell equation.

\section{Main Result} 
{\em For $k \geq 3$, the numbers that are both $k$-gonal and centered $k$-gonal are given by
$$\frac{k}{16(k - 2)}\left\{\frac{-2k^2 + 18k - 32}{k} + \left[k - 1 + \sqrt{k(k - 2)}\right]^{2i + 1} + \left[k - 1 - \sqrt{k(k - 2)}\right]^{2i + 1}\right\}$$
for $ i \geq 0$.}
\vskip 0.5cm
{\bf Proof.} Solving $P(n; k) = C(m; k)$ for $n$, we get
$$n = \frac{k - 4 + \sqrt{(4 - k)^2 -4(k - 2)(-km^2 + km - 2)}}{2(k - 2)}.$$
The expression under the radical can be rewritten as $k(k - 2)(2m - 1)^2 + 2k$ and must be a perfect square. 
Set 
\begin{equation}
k(k - 2)(2m - 1)^2 + 2k = a^2.
\end{equation}
This implies $k$ divides $a^2$. We consider two cases.
\vskip 0.5cm
{\bf Case 1:} If $a^2 = k^2b^2$, then $(k - 2)(2m - 1)^2 - kb^2 = -2$. This can be written as the Pell equation
\begin{equation}
\left[\frac{kb - (k - 2)A}{2}\right]^2 - k(k - 2)\left[\frac{(A - b)}{2}\right]^2 = 1,
\end{equation}
where $A = 2m - 1$.

If $k$ is odd, let 
$$Z = \frac{kb - (k - 2)A}{2}$$
and 
$$W = \frac{A - b}{2}.$$
Then (2) can be written as
 $$Z^2 - k(k - 2)W^2 = 1.$$
The ``smallest" solution is $(Z_0, W_0) = (k - 1, 1)$. So the general solution is given by (see [3], for instance)
$$Z = \frac{1}{2}\left\{\left[ k - 1 + \sqrt{k(k - 2)}\right]^i + \left[k - 1 - \sqrt{k(k - 2)}\right]^i\right\}$$
and
$$W = \frac{1}{2\sqrt{k(k - 2)}}\left\{\left[ k - 1 + \sqrt{k(k - 2)}\right]^i - \left[k - 1 - \sqrt{k(k - 2)}\right]^i\right\}.$$

Solving for $A$ and $b$, we get $A = Z + kW$ and $b = Z + (k - 2)W$. Letting $\alpha = k - 1 + \sqrt{k(k - 2)}$ and $\beta = k - 1 - \sqrt{k(k - 2)}$, we get the general terms $A_i$ and $b_i$ of the sequence of solutions to be

$$A_i  =  \frac{1}{2}\left(\alpha^i + \beta^i\right) + \frac{k}{2\sqrt{k(k - 2)}}\left( \alpha^i - \beta^i\right)$$
and 
$$b_i = \frac{1}{2}\left(\alpha^i + \beta^i\right) + \frac{k - 2}{2\sqrt{k(k - 2)}}\left(\alpha^i - \beta^i\right).$$
Since $A = 2m - 1$,
\begin{equation}
m_{(k, i)} = \frac{1}{4}\left[2 + \frac{k + \sqrt{k(k - 2)}}{\sqrt{k(k - 2)}}\alpha^i + \frac{-k + \sqrt{k(k - 2)}}{\sqrt{k(k - 2)}}\beta^i\right],
\end{equation}
where $m_{(k, i)}$ is the term that corresponds to $A_i$ for a given $k$. \Bigg(If $k$ is even, we let $Z = \frac{kb - (k - 2)A}{2}$ and $W = A - b$. Then
$$ Z^2 - \frac{k(k - 2)}{4}W^2 = 1.$$
the ``smallest" solution is $(Z_0, W_0) = (k - 1, 2)$ and so the general solution is given by
$$Z = \frac{1}{2}\left\{\left[ k - 1 + \frac{2\sqrt{k(k - 2)}}{2}\right]^i + \left[k - 1 - \frac{2\sqrt{k(k - 2)}}{2}\right]^i\right\}$$
and

$$W = \frac{1}{\sqrt{k(k - 2)}}\left\{\left[ k - 1 + \sqrt{k(k - 2)}\right]^i - \left[k - 1 - \sqrt{k(k - 2)}\right]^i\right\}.$$
Solving for $A$ and $b$, we get $A = Z + \frac{k}{2}W$ and $b = Z + \frac{k - 2}{2}W$. This gives the same values for $A_i$ and $b_i$ found in the odd case.\Bigg)

Expanding and  taking into account that $\alpha\beta = 1$ and $(\alpha + 1)(\beta + 1) = 2k$, we obtain
\begin{eqnarray*}
m_{(k, i)}^2 & = & \frac{1}{16} \Biggr\{4 + \frac{\left[k + \sqrt{k(k - 2)}\right]^2}{k(k - 2)}\alpha^{2i} + 
\frac{\left[-k + \sqrt{k(k - 2)}\right]^2}{k(k - 2)}\beta^{2i} + 4\frac{k + \sqrt{k(k - 2)}}{\sqrt{k(k - 2)}}\alpha^i \\
 & + & 4\frac{-k + \sqrt{k(k - 2)}}{\sqrt{k(k - 2)}}\beta^i - \frac{4}{k - 2} \Biggr\}.
\end{eqnarray*}
It follows, after simplifying, 
\[
C(m_{(k, i)}; k) = \frac{km_{(k, i)}^2 - km_{(k, i)} + 2}{2} = \frac{k}{16(k - 2)}\left\{\frac{-2k^2 + 18k - 32}{k} + \alpha^{2i + 1} + \beta^{2i + 1}\right\}.
\]
We now show that 
$$
n = \frac{k - 4 + a}{2(k - 2)},
$$
where $a$ is defined in (1), is actually an integer. In fact, substituting the value of $m_{(k, i)}$ found in (3) for $m$  and simplifying, we obtain
$$
a = a_i = \sqrt{k(k - 2)(2m - 1)^2 + 2k}
= \frac{\alpha^{i + 1}}{2} + \frac{\alpha^i}{2} + \frac{\beta^{i + 1}}{2} + \frac{\beta^i}{2},$$
where we have used the identitities $\left[k + \sqrt{k(k - 2)}\right]^2 = 2k\alpha$, $\left[k - \sqrt{k(k - 2)}\right]^2 = 2k\beta$, $k + \sqrt{k(k - 2)} = \alpha + 1$, and $k - \sqrt{k(k - 2)} = \beta + 1.$ It follows that 
\begin{eqnarray*}
n & = & \frac{k - 4 + a}{2(k - 2)}\\
& = & \frac{k - 4 + \frac{\alpha^{i + 1}}{2} + \frac{\alpha^i}{2} + \frac{\beta^{i + 1}}{2} + \frac{\beta^i}{2}}{2(k - 2)}.
\end{eqnarray*}
To show this is an integer, we notice first that for $i = 0$, 
$$n = \frac{2k - 8 + 2k - 2 + 2}{4(k - 2)} = 1.$$
Now
\begin{eqnarray*}
a_{i + 1} - a_i & = & \frac{\alpha^{i + 2} - \alpha^i + \beta^{i + 2} - \beta^i}{2}\\
& = & \frac{\alpha^i(\alpha^2 - 1) + \beta^i(\beta^2 - 1)}{2}\\
& = & \alpha^i\left[k(k - 2) + (k - 1)\sqrt{k(k - 2)}\right] + \beta^i\left[k(k - 2) - (k - 1)\sqrt{k(k - 2)}\right]\\
& = & \alpha^{i + 1}\sqrt{k(k - 2)} - \beta^{i + 1}\sqrt{k(k - 2)}\\
& = & \sqrt{k(k - 2)}\left(\alpha^{i + 1} - \beta^{i + 1}\right)\\
& = & \sqrt{k(k - 2)}\left(\alpha - \beta\right)\sum_{k=0}^i\alpha^{i-k}\beta^{k}\\
& = & \sqrt{k(k - 2)}\left(\alpha - \beta\right)h_i\\
& = & 2k(k - 2)h_i.
\end{eqnarray*}
Since $\alpha\beta = 1$, and $\alpha^i + \beta^i $ is an integer (see the value of $Z$ above),  $h_i$ is an integer and $a_{i+1} - a_i$ is a number divisible by $2k(k - 2)$. A routine induction argument on $i$ now shows that $\frac{k - 4 + a_i}{2(k - 2)}$ is an integer for all positive integers $i$. This completes the first case.

{\bf Case 2.} $a^2 = kb^2$, where $k$ is a perfect square. The equation reduces to

\begin{equation}
b^2 - (k - 2)A^2   = 2,
\end{equation}
where $A = \sqrt{k}(2m - 1)$. The solution to this equation, where $k - 2$ is not a perfect square, is standard [3]. We provide a sketch of the proof. First notice that the set of generating solutions contains $(\sqrt{k}, 1)$. To get all other solutions, we solve the Pell equation $b^2 - (k - 2)A^2 = 1$. If $(u, v)$ is a solution to the standard Pell equation, then $(\sqrt{k} + (k - 1)uv, \sqrt{k}v + u)$ is a solution to the same equation with the number 1 on the right hand side of the equation is  replaced with 2. Now a starightforward computation yields the same solutions we obtained in case 1. 
\vskip 0.5cm
\noindent

References

1. American Mathematical Monthly [1951, p. 568], problem E 954.

2. H. Eves, An Introduction to the History of Mathematics, revised edition, Holt-Winston.

3. K. H. Rosen, Elementary Number Theory and Its Applications, 5th edition, Pearson, 2005.

\end{document}